\documentclass[12pt]{article}
\usepackage{graphicx}
\usepackage{graphics}
\usepackage{amsthm}
\usepackage{amssymb, amsmath}
\title{Some integrals of the Dedekind $\eta$-function} 
\author{Mark W. Coffey\\
Department of Physics\\
Colorado School of Mines\\
Golden, CO  80401\\
mcoffey@mines.edu\\
Department of Mathematics\\
University of Colorado\\
Boulder, CO  80309\\
USA\\
mcoffey@colorado.edu\\}
\date{December 27, 2018}
\begin{document}
\maketitle
\baselineskip=25 pt
\begin{abstract}

Let $\eta$ be the weight $1/2$ Dedekind function.  A unification and generalization
of the integrals $\int_0^\infty f(x)\eta^n(ix)dx$, $n=1,3$, of Glasser \cite{glasser2009} is presented.  Simple integral inequalities as well as
some $n=2$, $4$, $6$, $8$, $9$, and $14$ examples are also given.  A prominent result is that
$$\int_0^\infty \eta^6 (ix)dx= \int_0^\infty x\eta^6 (ix)dx
={1 \over {8\pi}}\left({{\Gamma(1/4)} \over {\Gamma(3/4)}}\right)^2,$$
where $\Gamma$ is the Gamma function.
The integral $\int_0^1 x^{-1} \ln x ~\eta(ix)dx$ is evaluated in terms of a 
reducible difference of pairs of the first Stieltjes constant $\gamma_1(a)$.

\end{abstract}
 
\medskip
\baselineskip=15pt
\centerline{\bf Key words and phrases}
\medskip 

Dedekind function, integral, functional equation, Dirichlet $\beta$ function,
Gamma function, Gauss hypergeometric function, complete elliptic integral of the
first kind 

\bigskip
\noindent
{\bf 2010 AMS codes}
\newline{11F20, 30E20}

\baselineskip=25pt

\pagebreak
\centerline{\bf 1.  Introduction and preliminaries} 
\medskip

The Dedekind $\eta$ function is a modular form of weight $1/2$ and level $1$, defined on 
the upper half plane of $\mathbb{C}$, having general functional equation
$$\eta\left({{a\tau+b} \over {c\tau+d}}\right)=\epsilon(c\tau+d)^{1/2}\eta(\tau),$$
where $\epsilon^{24}=1$, $a$, $b$, $c$, and $d$ are integers, and $ad-bc=1$.  
Indeed, $\epsilon=\epsilon(a,b,c,d)$ may be made explicit (e.g., Theorem 3.4 of \cite{apostolmodfcns}).  If $c=0$, then simply $\epsilon(1,b,0,1)=\exp(i\pi b/12)$.
\footnote{Then the identity transformation of the modular group $\Gamma$ with $b=0$ is
obvious.}
If $c>0$, $\epsilon$ may be expressed in terms of a Dedekind sum. 
\footnote{However, there is no obvious doubling formula relating $\eta(2\tau)$ to
$\eta(\tau)$.}

It follows that
$$\eta(\tau+1)=e^{i\pi/12}\eta(\tau), ~~~~~~\eta\left(-{1 \over \tau}\right)
=\sqrt{-i\tau}\eta(\tau). \eqno(1.1)$$
The latter relation will specifically be useful later. 
It follows that there are special values, for instance, for $\eta(ji)$ and
$\eta(i/j)$, $j \in \mathbb{N}^+$.  Many of these values, including $\eta(ji)$,
involve the factor $\Gamma(1/4)$, where $\Gamma$ is the Gamma function.

With $q=\exp(-2\pi x)$, $\eta$ has the product development
$$\eta(ix)=q^{1/24}\prod_{n=1}^\infty (1-q^n).$$
This function has the series representation \cite{andrews}
$$\eta(ix)={2 \over \sqrt{3}}\sum_{n=0}^\infty \cos[(2n+1)\pi/6]q^{(2n+1)^2/24}
=q^{1/24}\sum_{n=-\infty}^\infty (-1)^n q^{n(3n-1)/2},
\eqno(1.2)$$
and via Jacobi's triple identity,
$$\eta^3(ix)=\sum_{n=0}^\infty (-1)^n(2n+1)q^{(2n+1)^2/8}
=q^{1/8}\sum_{n=0}^\infty (-1)^n(2n+1)q^{n(n+1)/2}. \eqno(1.3)$$
We may note that in (1.2), $\sqrt{3}\cos[(2n+1)\pi/6]/2$ is similar to a 
real-valued Dirichlet character $\chi_{12}(n)$ modulo $12$.
The $\eta$ function is related to several theta series.
Example sources of further information on $\eta$ are Ch.\ 8.8 of \cite{rroy2017}
and Ch.\ 3 of \cite{apostolmodfcns}.

The Dedekind eta function is of interest in number theory, combinatorics, and
conformal field theory (e.g., \cite{andrews,fuchs,kac}).  An example is that by applying the Weyl-Kac character formula and the denominator identities, it is possible to find numerous combinatorial identities.  We mention such identities further in section 2. 

In order to more easily connect with other areas, we mention $q$ series and $q$ product
notation for the $\eta$ function.  Let $(\alpha)_m=(\alpha;q)_m=\prod_{i=0}^{m-1} (1-\alpha
q^i)$, with $|q|<1$ and $m \geq 0$ an integer, so that $\eta(\tau)=q^{1/24}(q;q)_\infty$,
where $\tau=\ln q/(2\pi i)$.  Then the generating function for the number of partitions of
$n$, $p(n)$, is given by
$$\sum_{n=0}^\infty p(n)q^n={1 \over {(q;q)_\infty}}={q^{1/24} \over {\eta(\tau)}}.$$
From the Durfee square there is the identity 
$$\sum_{n=0}^\infty {q^{n^2} \over {(q;q)_n^2}}={1 \over {(q;q)_\infty}}={q^{1/24} \over
{\eta(\tau)}}=\prod_{n=1}^\infty (1-q^n)^{-1},$$
and this extends to
$$\sum_{n=0}^\infty {z^nq^{n^2} \over {(q;q)_n (zq;q)_n}}={1 \over {(zq;q)_\infty}}.$$

In \cite{glasser2009}, Glasser used Laplace and Fourier transforms to give some integrals $\int_0^\infty f(x)\eta^n(ix)dx$, for $n=1$ and $n=3$.  Other such
integrals were presented in the Appendix of \cite{glasser2009} without further
elaboration.  In this paper, we especially unify and generalize many of the 
entries of that Appendix.

Accordingly, we also introduce the Dirichlet $\beta$ function
$$\beta(s)=\sum_{n=0}^\infty {{(-1)^n} \over {(2n+1)^s}}
=4^{-s}\left[\zeta\left(s,{1\over 4}\right)-\zeta\left(s,{3\over 4}\right)\right],$$
where $\zeta(s,a)$ is the Hurwitz zeta function.  The above series holds
for Re $s>0$, but $\beta$ is extended to $\mathbb{C}$ through analytic
continuation, with the following functional equation holding:
$$\beta(1-s)=\left({\pi \over 2}\right)^{-s} \sin\left({{\pi s} \over 2}\right)
\Gamma(s)\beta(s).$$
The famous and ubiquitous value $\beta(2)=G$ is the Catalan constant, such that
\cite{summatory}
$$G={\pi^2 \over 8}\prod_{\overset{p \equiv 3}{mod ~4}} {{p^2-1} \over {p^2+1}}=
{\pi^2 \over {12}}\prod_{\overset{p \equiv 1}{mod ~4}} {{p^2+1} \over {p^2-1}}
={1 \over 2}\int_0^\infty {x \over {\cosh x}}dx$$
$$={1 \over {16\pi^2}}\sum_{k=0}^\infty {1 \over {k!}}\left[\gamma_k\left({1 \over 4}\right)-\gamma_k\left({3 \over 4}\right)\right]={1 \over {16}}\left[\psi'
\left({1 \over 4}\right)-\psi'\left({3 \over 4}\right)\right],$$
and $\gamma_k(a)$ are the Stieltjes constants \cite{summatory,coffeyrmjm2011,coffeyrama} to be recollected in section 5.
The products in this equation are taken over prime numbers $p$.

We recall the polygamma functions $\psi^{(j)}(z)=d^{j+1}\ln \Gamma/dz^{j+1}$.
In particular for the digamma and trigamma functions at positive integer argument
$\ell$,
$$\psi(\ell)=-\gamma+\sum_{k=1}^{\ell-1} {1 \over k}, ~~~~~~
 \psi'(\ell)={\pi^2 \over 6}-\sum_{k=1}^{\ell-1} {1 \over k^2},$$
where $\gamma=-\psi(1)$ is the Euler constant.  More generally, the polygamma
functions at positive integer argument evaluate in terms of values of the
Riemann zeta function $\zeta(s)=\zeta(s,1)$ and generalized harmonic numbers
$H_n^{(r)}=\sum_{k=1}^n 1/k^r$:
$$\psi^{(j)}(n+1)=(-1)^{j+1}j!\left[\zeta(j+1)-H_n^{(j+1)}\right],$$
wherein $\psi^{(j)}(1)=(-1)^{j+1}j!\zeta(j+1)$.
In addition, $\psi^{(j)}(1/2)$ may be determined in terms of the zeta function 
values $\zeta(j+1)$.  The polygamma functions have several functional equations
and many integral representations, including for $j>0$
$$\psi^{(j)}(z)=\int_0^1 {{t^{z-1}\ln^j t} \over {t-1}}dt.$$

Reference \cite{patkowski} has considered integrals 
$$\int_0^\infty \eta^3\left({{i4x} \over \pi}\right)e^{-b^2 x}t(x)dx,
~~~~\int_0^\infty \eta^6\left({{i4x} \over \pi}\right)t(x)dx,$$
wherein $t(x)=\sin cx$ or $\cos cx$.  The method there is the application of
Poisson summation.  Below we mention a connection with our results in the
very special case of $c=0$ in Theorem 1.2 in \cite{patkowski}.  While, for
example, (3.4) there is referred to as Poisson summation, the trapezoidal
sum $T=\sum'_{a\leq n \leq b} f(n)$ readily follows from Euler summation 
according to the Fourier series for the periodized first and second Bernoulli
polynomials $B_1(x-[x])$ and $B_2(x-[x])$, where $[\mbox{.}]$ is the greatest 
integer function.

No integrals with $\eta^2(ix)$, $\eta^4(ix)$, $\eta^6(ix)$, or $\eta^8(ix)$ in the integrand
appear in \cite{glasser2009}.  
We provide such example results in the next section.  Then we give
another major result which unifies, depending upon the point of view,
at least three of the entries of the Appendix of \cite{glasser2009}.  
In a later section we show several equivalences of pairs or of triples, even of
quintuples, of entries of that Appendix.

This paper often uses the interchange of infinite series and integral.  The various 
series involve functions $\sim e^{-k n^2 x}$ with a prefactor $O(n)$ and fixed $k>0$ 
as $n \to \infty$.
As such, the interchange is justified on the basis of Levi's theorem for series of
Lebesgue-integrable functions.  

\medskip
{\bf 2.1 Integrals $\int_0^\infty \eta^k(ix)dx$, $k=2, 4, 6, 8, 9$}
\medskip

The evaluations $\int_0^\infty \eta(ix)dx=2\pi/\sqrt{3}$ and $\int_0^\infty \eta^3(ix)dx=1$ are known.  We first provide an expression for $\int_0^\infty \eta^2(ix) dx$ in terms of a single-index infinite trigonometric series.
Then we discuss and provide specific expressions for  
$\int_0^\infty x^j\eta^6(ix)dx$.

{\bf Proposition 1}.  Let
$$C \equiv \sum_{m=0}^\infty {{\cos[(2m+1)\pi/6]} \over {2m+1}}\left\{\cot[(-5+i(2m+1))\pi/12] + \tan[(-5+i(2m+1))\pi/12]\right\}.$$
Then
$$\int_0^\infty \eta^2(ix) dx ={1 \over \sqrt{3}}\mbox{Im} ~C.$$

From (1.2) we have
$$\int_0^\infty \eta^2(ix) dx={4 \over 3}\sum_{n,m=0}^\infty \cos\left[{{(2n+1)\pi} \over 6}\right]
\cos\left[{{(2m+1)\pi} \over 6}\right]\int_0^\infty q^{[(2n+1)^2+(2m+1)^2]/24}dx$$
$$={8 \over \pi}\sum_{n,m=0}^\infty {{\cos[(2n+1)\pi/6] \cos[(2m+1)\pi/6]} \over
{2n(n+1)+2m(m+1)+1}}.$$
The sum over $n$ or $m$ may be performed in terms of differences of pairs of
digamma functions $\psi=\Gamma'/\Gamma$.  With $\psi(z)-\psi(1-z)=-\pi \cot \pi z$, there results
$$\int_0^\infty \eta^2(ix) dx ={1 \over \sqrt{3}}{1 \over i}(C-C^*),$$
hence the evaluation. \qed

{\bf Proposition 2}. (a) Let
$$C_2 \equiv \sum_{m=0}^\infty (-1)^m (2m+1)\mbox{sech}\left({{(2m+1)\pi} \over 2}
\right)$$
$$={1 \over {8\pi}}\left({{\Gamma(1/4)} \over {\Gamma(3/4)}}\right)^2.$$
Then
$$\int_0^\infty \eta^6(ix) dx = C_2.$$

(b)
$$\int_0^\infty x\eta^6(ix) dx =\int_0^\infty \eta^6(ix) dx = C_2.$$

(c) Let
$$C_3={1 \over {2\pi}}\sum_{m=0}^\infty {{(-1)^m\mbox{sech}^3[(m+1/2)\pi)((2m+1)\pi (-3+\cosh(2m+1)\pi)+2\cosh((2m+1)\pi)]} \over {(2m+1)^2}}.$$

Then
$$\int_0^\infty x^2\eta^6(ix) dx = C_3.$$

(d) Let $a,b>0$.  Then the integral
$$\int_0^\infty \eta^3(iax)\eta^3(ibx)dx$$
may be expressed in closed form.

(e) Let
$$C_4\equiv {2 \over \sqrt{3}}\sum_{n=0}^\infty \cos\left[(2n+1){\pi \over 6}
\right]\mbox{sech}\left[{{(2n+1)\pi} \over {2\sqrt{3}}}\right]$$
$$=\sum_{\overset{n \geq 0}{n \equiv 0,2,3,5,6,8,9,11 \mod 12}} (\pm 1) \mbox{sech}\left[{{(2n+1)\pi} \over {2\sqrt{3}}}\right]$$
$$=\sum_{n=-\infty}^\infty (-1)^n \mbox{sech}\left[{{(6n+1)\pi} \over {2\sqrt{3}}}
\right].$$
Then
$$\int_0^\infty \eta^4(ix)dx=C_4.$$

(f) We have the identity
$$5\int_0^\infty \eta^8(iy)dy=e^{-i\pi/3}\int_0^\infty \eta^8\left(ix+{1 \over 2}\right)dx.$$

(a) The series follows from the use of (1.3), so that
$$\int_0^\infty \eta^6(ix) dx ={2 \over \pi}\sum_{n,m=0}^\infty {{(-1)^{n+m}(2n+1)(2m+1)} \over {2n(n+1)+2m(m+1)+1}}.$$ 

For the closed form evaluation, we put $y=\pi$ in Entry 16(i) of Ch.\ 17 of \cite{ramapart3}.
Then
$$C_2={z^2 \over 2}\sqrt{x(1-x)},$$
where $x=1/2$ and 
$$z={}_2F_1\left({1 \over 2},{1 \over 2};1;{1 \over 2}\right)
={{\Gamma(1/4)} \over {\sqrt{2\pi} \Gamma(3/4)}},$$
and $_2F_1$ is the Gauss hypergeometric function \cite{a2r,nbs,grad}. 

(b) This result is equivalent to the summation identity
$$\sum_{n,m=0}^\infty {{(-1)^{n+m}(2n+1)(2m+1)} \over {[2n(n+1)+2m(m+1)+1]}}
={2 \over \pi}\sum_{n,m=0}^\infty {{(-1)^{n+m}(2n+1)(2m+1)} \over {[2n(n+1)+2m(m+1)+1]^2}},$$
or
$$\sum_{m=0}^\infty (-1)^m (2m+1)\mbox{sech}\left({{(2m+1)\pi} \over 2}\right)
={\pi \over 8}\sum_{m=0}^\infty (-1)^m \left[\csc^2{1 \over 4}(1+i(2m+1)\pi)
-\sec^2{1 \over 4}(1+i(2m+1)\pi)\right].$$

(c) Generally for integer $j\geq 0$,
$$\int_0^\infty x^j\eta^6(ix) dx ={{2^{j+1}j!} \over \pi^{j+1}}\sum_{n,m=0}^\infty {{(-1)^{n+m}(2n+1)(2m+1)} \over {[2n(n+1)+2m(m+1)+1]^{j+1}}}.$$ 

When summed over $n$, pairs of differences of polygamma functions $\psi^{(j)}$,
and of lower orders appear.  The arguments of these functions are 
$${1 \over 4}\pm {i \over 4}(2m+1) ~~~~\mbox{and} ~~~~ {3 \over 4}\pm {i \over 4}(2m+1).$$
As such, these difference pairs may be reduced.

In the case of $j=2$,
$$\int_0^\infty x^2\eta^6(ix) dx ={{16} \over \pi^3}\sum_{n,m=0}^\infty {{(-1)^{n+m}(2n+1)(2m+1)} \over {[2n(n+1)+2m(m+1)+1]^3}}.$$ 
When summed over $n$, there are two pairs of $\psi'$ and six pairs of $\psi''$
at the aforementioned arguments.  Upon simplification, $C_3$ results.  

(d)  We have
$$\int_0^\infty \eta^3(iax)\eta^3(ibx)dx=\sum_{n,m=0}^\infty (-1)^{n+m} (2n+1)
(2m+1) \int_0^\infty q^{(2n+1)^2/8}(ax)q^{(2m+1)^2/8}(bx)dx$$
$$={4 \over \pi}\sum_{n,m=0}^\infty {{(-1)^{n+m} (2n+1) (2m+1)} \over {a(2n+1)^2+b(2m+1)^2}}$$
$$={1 \over a}\sum_{m=0}^\infty (-1)^m (2m+1)\mbox{sech}\left(\sqrt{{b \over a}}
(2m+1){\pi \over 2}\right).$$
We now use entry 16(i) of Ch.\ 17 of \cite{ramapart3} with 
$y=\sqrt{{b \over a}}\pi$.  The result of the integration is
$${1 \over {2a}}z^2\sqrt{x(1-x)},$$
where $x$ is found from $y$ via
$$y=\pi {{{}_2F_1(1/2,1/2;1;1-x)} \over {{}_2F_1(1/2,1/2;1,x)}}.$$
Then $k=\sqrt{x}$, $k'=\sqrt{1-x}$, $K=K(k)$, $K'=K(k')$, $y=\pi K'/K$, and 
$z={2 \over \pi}K$.

(e) Method 1.  From (1.2) and (1.3) we obtain
$$\int_0^\infty \eta^4(ix)dx={2 \over \sqrt{3}}{3 \over \pi}\sum_{m,n=0}^\infty
\cos\left[(2n+1){\pi \over 6} \right]{{(-1)^m(2m+1)} \over {[3m^2+n^2+3m+n+1]}}
=C_4.$$

Method 2.  We may use the Laplace transform $\int_0^\infty e^{-xy}\eta^3(ix)dx$
((14) or (A.1) in \cite{glasser2009}) and the unilateral and bilateral series
in (1.2).  For instance, putting $y\to (2n+1)^2\pi/12$, and summing over $n$ with the
$(2/\sqrt{3})\cos[(2n+1)\pi/6]$ coefficient,
$${2 \over \sqrt{3}}\sum_{n=0}^\infty \cos\left[(2n+1){\pi \over 6}\right] \int_0^\infty e^{-(2n+1)^2\pi x/12} \eta^3(ix)dx=\int_0^\infty \eta^4(ix)dx$$
$$={2 \over \sqrt{3}}\sum_{n=0}^\infty \mbox{sech}\left[{{(2n+1)\pi} \over 
{2\sqrt{3}}}\right]\cos\left[(2n+1){\pi \over 6}\right].$$

Method 3.  We may put $c \to 0$ in (1.4) of Theorem 1.1 of \cite{patkowski}, in 
which case $A=b$ and $B=0$ there, giving the Laplace transform
$$\int_0^\infty \eta^3(iy)e^{-b^2\pi y/4}dy={1 \over {\cosh \pi b/2}}.$$
Then we use the bilateral series in (1.2) with $b^2=4n(3n-1)+1/3$.  Then summing over
$n$ we obtain
$$\int_0^\infty \eta^4(iy)dy=\sum_{n=-\infty}^\infty {{(-1)^n} \over {\cosh\left[{\pi \over 2}\sqrt{4n(3n-1)+1/3}\right]}}$$
$$=\sum_{n=-\infty}^\infty {{(-1)^n} \over {\cosh\left[{\pi \over {2\sqrt{3}}}(6n-1)\right]
}}. $$

(f) This follows from application of the identity coming from either the product or
series developments of $\eta$,
$$\eta^8(\tau)+16\eta^8(4\tau)=e^{-i\pi/3}\eta^8\left(\tau+{1 \over 2}\right).$$

\qed  

{\it Remarks}.  An alternative series for $C_2$ is
$$C_2={1 \over \pi}\sum_{m=0}^\infty (-1)^m (2m+1)\Gamma\left[{1 \over 2}+{i \over 2}
(2m+1)\right]\Gamma\left[{1 \over 2}-{i \over 2} (2m+1)\right].$$

The value of $z$ above in part (a) comes from the complete elliptic integral of the first kind $K=K(1/\sqrt{2})$.  When $y=\pi$ as above, the modulus $k=1/\sqrt{2}=k'$,
the complementary modulus.

The summation of part (a) agrees with the degenerate case of $c=0$, $A(n,c=0)=n$,
and $B=0$ in (1.7) of \cite{patkowski}.  Then in the notation of that reference
with $\chi(n)=\chi_4(n)$,
$$\int_0^\infty \eta^6(iy)dy=\sum_{n=1}^\infty {{\chi(n)n\cosh(\pi n/2)} \over
{\cosh^2(\pi n/2)}}$$
$$=\sum_{n=1}^\infty {{\chi(n)n} \over {\cosh(\pi n/2)}}
=\sum_{k=0}^\infty {{(2k+1)\chi(2k+1)} \over {\cosh \pi(2k+1)/2}}$$
$$=\sum_{k=0}^\infty {{(2k+1)(-1)^k} \over {\cosh \pi(2k+1)/2}}.$$
Reference \cite{patkowski} mentions that the right sides of its (1.6) and (1.7)
\footnote{Otherwise these equations are referred to as (1.4) and (1.5).}
are not expressible in terms of elementary functions.  But this is entirely
expected, and does not rule out the possibility of closed form evaluation in
certain instances in terms of special functions as parts (a) and (d) indicate.

Part (d) is a partial answer to the challenge posed at the end of \cite{patkowski}.
In addition, based upon the Weyl-Kac character formula, we have identities such as
(\cite{fuchs}, p.\ 137)
$$\eta(8\tau)\eta(16\tau)=\sum_{\overset{m,n\in \mathbb{Z}}{m \leq |3n|}} (-1)^m
q^{(2m+1)^2-32n^2},$$
$$\eta(12\tau)\eta(12\tau)=\sum_{\overset{m,n\in \mathbb{Z}}{m \leq |2n|}} 
(-1)^{m+n} q^{(3m+1)^2-(6n+1)^2},$$
and
$$\eta(24\tau)\eta(96\tau)=\sum_{\overset{m,n\in \mathbb{Z}}{2m \leq n \leq 0}} 
(-1)^{n(n+1)/2} q^{8(3m+1)^2-3(2n+1)^2}(1-q^{24(2m+1)}).$$

These lead to summation expressions and integral inequalities. 
For instance, we have
$${1 \over {2\pi}}-{1 \over {8\pi}}\left(8-\sqrt{2}\pi\cot\left({\pi \over {4\sqrt{2}}}\right)\right)<\int_0^\infty \eta(8ix)\eta(16ix)dx< {1 \over {2\pi}}.$$
This follows by separating the $m=n=0$ term, so that
$$\int_0^\infty \eta(8ix)\eta(16ix)dx={1 \over 8}\int_0^\infty \eta(iy)\eta(2iy)dy
=\int_0^\infty q(x)dx+2\sum_{n=1}^\infty
\sum_{m=-3n}^{3n}(-1)^m \int_0^\infty q^{(2m+1)^2-32n^2}dx$$
$$={1 \over {2\pi}}+{2 \over \pi} \sum_{n=1}^\infty \sum_{m=-3n}^{3n} {{(-1)^m} 
\over {[(2m+1)^2-32n^2]}}.$$
The left inequality follows by applying the result of partial fractions so that,
from the $m=0$ term,
$$\sum_{n=1}^\infty {1 \over {32n^2-1}}={1 \over {16}}\left(8-\sqrt{2}\pi\cot\left(
{\pi \over {4\sqrt{2}}}\right)\right).$$

Similarly,
$$\int_0^\infty \eta^2(12ix) dx < {1 \over {2\pi}}.$$
In fact, for the latter integral, from Proposition 1 we have
$$\int_0^\infty \eta^2(12ix) dx = {1 \over {12}}\int_0^\infty \eta^2(ix) dx
={1 \over {12\sqrt{3}}}\mbox{Im} ~C.$$


{\bf Proposition 3}.  Let
$$C_9 \equiv \sum_{m,n=0}^\infty (-1)^{m+n}(2m+1)(2n+1)\mbox{sech}\left(\sqrt{2m
(m+1)+2n(n+1)+1}{\pi \over \sqrt{2}}\right).$$
Then
$$\int_0^\infty \eta^9(ix)dx=C_4.$$
We note that $C_9$ may be reducible to a one-dimensional sum via the use of
entry 16(i) of Ch.\ 17 \cite{ramapart3}.  This is in reference to the use of the values $K$ and $K'$ of the complete elliptic integral of the first kind.

We have
$$\int_0^\infty \eta^9(ix)dx={4 \over \pi}\sum_{\ell,m,n=0}^\infty {{(-1)^{\ell+m+n}
(2\ell+1)(2m+1)(2n+1)} \over {[4\ell(\ell+1)+4m(m+1)+4n(n+1)+3]}}.$$
After summing over $\ell$ and simplifying, we obtain $C_9$. \qed

Omitting the proof, we mention the inequalities contained in the following.
{\newline \bf Proposition 4}. (a) For $j > 1$, $\int_0^\infty \eta^{j-1}(ix)dx
\geq \int_0^\infty \eta^j(ix)dx$ and (b) $\int_0^\infty \eta^j(ix)dx \leq {{12}
\over {\pi j}}$.  (c) More generally than (a), for $j>1$ and $n \geq 0$,
$\int_0^\infty x^n\eta^{j-1}(ix)dx \geq \int_0^\infty x^n\eta^j(ix)dx$. 

\medskip
{\bf 2.2 The value ${}_2F_1(1/2,1/2;1,1/2)$}
\medskip

We discuss this factor which occurs in $C_2$ (cf.\ \cite{a2r}, pp.\ 126-128, 136).
Through transformation formulae for ${}_2F_1$, only herein very briefly mentioned
(e.g., \cite{grad}, p.\ 1043), values $_2F_1(z)$ relate to those of ${}_2F_1[z/
(z-1)]$.  Therefore, $_2F_1(-1)$ is proportional to $_2F_1(1/2)$, and in the
subject case, Kummer's identity applies,
$${}_2F_1(a,b;a-b+1;-1)={{\Gamma(a-b+1)\Gamma(a/2+1)} \over {\Gamma(a+1)\Gamma(
a/2-b+1)}}.$$

We also have the following elegant specialization of $_2F_1(a,b;2a;1-x)$:
$${}_2F_1\left({1 \over 2},{1 \over 2};1;1-x\right)=\left({2 \over {1+x}}\right)^{1/2}
{}_2F_1\left[{1 \over 4},{3 \over 4};1;\left({{1-x} \over {1+x}}\right)^2\right].$$
There results several identities, including with the complete elliptic integral
of the first kind $K$:
$${}_2F_1\left({1 \over 2},{1 \over 2};1,{1 \over 2}\right)={2 \over \sqrt{3}} {}_2F_1\left({1 \over 4},{3 \over 4}; 1;{1 \over 9}\right)$$
$$={\sqrt{\pi} \over {\Gamma^2(3/4)}}={2 \over \pi} K\left({1 \over 2}\right).$$
There are several other quadratic transformations of the function ${}_2F_1$, but
we refrain from further elaboration.

\medskip
{\bf 3.  Generalization of A.2 and A.4 of \cite{glasser2009}}
\medskip

The following generalizes the two subject entries, with several corollaries,
including A.3.
{\bf Proposition 5}.  Let Re $p>0$.  Then
$$\int_0^\infty {{\eta^3(ix)dx} \over {(x+a)^p}}
=\int_0^\infty {{y^{p-1/2}\eta^3(iy)dy} \over {(1+ay)^p}}$$
$$={2 \over {\pi^p(p-1)!}}
\int_0^\infty {{x^{2p-1} e^{-ax^2/\pi}} \over {\cosh x}}dx.$$

In the Appendix A of \cite{glasser2009}, A.2 is the special case of $p=1$
and A.4 is the special case of $p=1/2$, when $(-1/2)!=\Gamma(1/2)=\sqrt{\pi}$.

{\bf Corollary 1}.  
$$\int_0^\infty x^{-p}\eta^3(ix)dx={4 \over \pi^p}{{\Gamma(2p)} \over {\Gamma(p)}}
\beta(2p).$$

This follows as the $a \to 0$ case of Proposition 4,
$$\int_0^\infty x^{-p}\eta^3(ix)dx={2 \over {\pi^p(p-1)!}}
\int_0^\infty {x^{2p-1} \over {\cosh x}}dx.$$

{\bf Corollary 2}.  Noting that
$$\Gamma(s)={1 \over s}-\gamma + O(s), ~~~~~\mbox{as} ~~ s \to 0,$$
$${{\Gamma(2p)} \over {\Gamma(p)}}=1/2+O(p), ~~\mbox{and} ~~ \beta(p)=1/2+O(p) ~~\mbox{as} ~~
p \to 0,$$
the $p \to 0$ limit of Corollary 1 gives
$$\int_0^\infty \eta^3(ix)dx=1,$$
as stated in A.14 of \cite{glasser2009}.

The importance of the variable $p$ in Proposition 4 includes that we can then
differentiate and/or integrate with respect to it, and/or sum over it.  An
example of the latter is the following.
\newline{\bf Corollary 3}.  Let $f(p)=\cos kp$ and $\sin kp$ and consider
$$\int_0^\infty \sum_{p=1}^\infty {{f(p) \eta^3(ix)} \over {(x+a)^p}}dx.$$
There results:
$$\int_0^\infty {{[(a+x)\cos k-1]} \over {1+(a+x)^2-2(a+x)\cos k}}\eta^3(ix)dx$$
$$={2 \over \pi}\int_0^\infty {{xe^{-ax^2/pi}} \over {\cosh x}}\cos(k+x^2 \sin k/\pi)
[\cosh(x^2 \cos k/\pi)+\sinh(x^2\cos k/\pi]dx$$
and
$${i \over 2}\int_0^\infty {{(e^{2ik}-1)(a+x)} \over {(e^{ik}-x-a)((a+x)e^{ik}-1)}}\eta^3(ix)dx$$
$$={2 \over \pi}\int_0^\infty {{xe^{-ax^2/\pi}} \over {\cosh x}}\sin(k+x^2 \sin k/\pi) [\cosh(x^2 \cos k/\pi)+\sinh(x^2\cos k/\pi]dx$$

In regard to Proposition 5, we briefly mention the use of binomial expansion
along with the use of the corrected form of A.15, the latter which is discussed
in the following section,
$$\int_0^\infty {{\eta^3(ix)dx} \over {(x+a)^p}}=\sum_{\ell=0}^\infty {{-p}
\choose \ell} a^{-p-\ell}\int_0^\infty \eta^3(ix)x^\ell dx$$
$$=\sum_{\ell=0}^\infty {{-p} \choose \ell} a^{-p-\ell}{{4^{\ell+1}\ell!} \over
\pi^{\ell+1}}\beta(2\ell+1).$$

There are then many routes to Proposition 5 as the integral representation
$$\zeta(s,a)={1 \over {\Gamma(s)}}\int_0^\infty {{x^{s-1}e^{-(a-1)x}} \over
{e^x-1}}dx,$$
for Re $s>1$ and Re $a>0$, gives for Re $s>0$
$$\beta(s)={1 \over {\Gamma(s)}}\int_0^\infty {{y^{s-1}e^y} \over {e^{2y}+1}}dy
={1 \over {2\Gamma(s)}}\int_0^\infty {{y^{s-1}dy} \over {\cosh y}}.$$ 

The second equality in the Proposition follows from the change of variable
$y=1/x$ and the use of (1.1).
\qed

We may then put $a=j$ and sum on $j$ from $0$ to $\infty$ in Proposition 5
to have
{\newline \bf Corollary 4}.
$$\int_0^\infty \eta^3(ix)\zeta(p,x)dx=\int_0^\infty y^{-1/2}\eta^3(iy)\zeta(p,1/y)
dy$$
$$={2 \over {\pi^p (p-1)!}}\int_0^\infty {{x^{2p-1} dx} \over {\cosh x (1-e^{-x^2 /\pi})}}.$$
In particular, this includes the cases for $n \geq 0$,
$\zeta(-n,x)=-B_{n+1}(x)/(n+1)$, where $B_n(x)$ is the $n$th Bernoulli polynomial.

\medskip
{\bf 4. Integral equivalences}
\medskip


{\bf Proposition 6}.  Integrals A.1, A.5, A.6, A.15, and A.3 of \cite{glasser2009} are equivalent.  Here
$$\int_0^\infty e^{-xy}\eta^3(ix)dx=\mbox{sech}\sqrt{\pi y}, \eqno(A.1)$$
being also (14) in the text of \cite{glasser2009}, 
$$\int_0^\infty x^{-1/2}e^{-a/x}\eta^3(ix)dx=\mbox{sech} \sqrt{\pi a}, \eqno(A.5)$$
and
$$\int_0^\infty e^{-xy}\eta^3(ix){{dx} \over x}={2 \over \pi}\int_{\sqrt{\pi y}}^\infty x\mbox{sech}x ~dx.  \eqno(A.6)$$

We have $\int_y^\infty e^{-xy'}dy'=e^{-xy}/x$ for Re $x>0$.  Accordingly
integrating both sides of (A.1) and making a simple substitution on the right
side shows the equivalence of (A.1) and (A.6).  
We next use the substitution $v=1/x$ in A.5 to write
$$\int_0^\infty x^{-1/2}e^{-a/x}\eta^3(ix)dx=\int_0^\infty v^{-3/2}e^{-av}
\eta^3\left({i \over v}\right)dv$$
$$=\int_0^\infty e^{-av}\eta^3(iv)dv=\mbox{sech} \sqrt{\pi a}.$$
In the second step, we employed the second functional equation on the right 
side of (1.1).  
The equivalence of A.3 and A.15 is shown in Proposition 8.

For the equivalence of A.1 and A.15 we use the generating function
$${{ze^z} \over {e^{2z}+1}}=\mbox{sech}z=\sum_{n=0}^\infty E_n {z^n \over {n!}},
~~~~~~|z|<\pi/2,$$
where $E_n$ are the Euler numbers.  The first nonzero few of these numbers are
$E_0=1$, $E_2=-1$, $E_4=5$, and $E_6=-61$.  Then by expanding both sides of A.1
in powers of $y$ and using that $E_{2n+1}=0$, we find that
$$\int_0^\infty x^n \eta^3(ix)dx=(-1)^n {{n!} \over {(2n)!}}\pi^n E_{2n}
={{4^{n+1}n!} \over \pi^{n+1}}\beta(2n+1).$$
Therefore, the equivalence of all five integrals is shown.  \qed

{\bf Proposition 7}.  Integrals A.8 and A.11 of \cite{glasser2009} are
equivalent.  Here
$$\int_0^\infty x^{-1/2}\cos\left({a \over x}\right)\eta^3(ix)dx=2{{\cos\sqrt{\pi
a/2}\cosh\sqrt{\pi a/2}} \over {\cos \sqrt{\pi a}+\cosh \sqrt{2\pi a}}}, 
\eqno(A.8)$$
and
$$\int_0^\infty \cos(xy)\eta^3(ix)dx={{\cosh\sqrt{\pi y/2}\cos\sqrt{\pi y/2}} \over
{\sinh^2 \sqrt{\pi y/2}+\cos^2 \sqrt{\pi y/2}}}.  \eqno(A.11)$$

We make the change of variable $v=1/x$ in A.8 and use (1.1) so that
$$\int_0^\infty x^{-1/2}\cos\left({a \over x}\right)\eta^3(ix)dx
=\int_0^\infty \cos (av) \eta^3(iv)dv.$$
Lastly there is the trigonometric identity $\sinh^2 x+\cos^2 x=(1/2)(\cosh 2x+
\cos 2x)$. \qed

{\bf Proposition 8}.  Integrals A.3 and A.15 (corrected) of \cite{glasser2009} are
equivalent.  Here, for $\nu>0$,
$$\int_0^\infty x^{-\nu}\eta^3(ix)dx={4 \over \pi^\nu}{{\Gamma(2\nu)} \over
{\Gamma(\nu)}}\beta(2\nu), \eqno(A.3)$$
and
$$\int_0^\infty x^n \eta^3(ix)dx={{4^{n+1}n!} \over \pi^{n+1}}\beta(2n+1).
\eqno(A.15)$$

Again the second functional equation of (1.1) is used, so that
$$\int_0^\infty x^{-\nu}\eta^3(ix)dx=\int_0^\infty v^{\nu-2}\eta^3\left({i \over v}
\right)dv =\int_0^\infty v^{\nu-1/2}\eta^3(iv)dv.$$
Then we put $n=\nu-1/2$ and apply the duplication formula for the Gamma 
function,
$$\Gamma(2n+1)={2^{2n} \over \sqrt{\pi}}\Gamma\left(n+{1 \over 2}\right) \Gamma(n+1).$$
\qed

{\it Remark}.  Similarly with the aid of (1.1) we may transform the left side
of A.10 to
$$\int_0^\infty x^{-1/2}e^{a/x} \mbox{erfc}\left(\sqrt{{a \over x}}\right)
\eta^3(ix)dx=\int_0^\infty v^{-3/2}e^{av}\mbox{erfc}(\sqrt{av})\eta^3\left({i \over
v}\right)dv$$
$$=\int_0^\infty e^{av}\mbox{erfc}(\sqrt{av})\eta^3(iv)dv,$$
where erfc is the complementary error function.


{\bf Proposition 9}. Entries A.11 and A.12 may obtained from entry A.15 (as
corrected).

We demonstrate this briefly obtaining an equivalent form of A.11 from A.15:
$$\int_0^\infty \cos (xy) \eta^3(ix)dx=\sum_{j=0}^\infty y^{2j} {{(-1)^j} \over {(2j)!}} \int_0^\infty x^{2j}\eta^3(ix)dx$$
$$=\sum_{j=0}^\infty (-1)^j y^{2j}\left({4 \over \pi}\right)^{2j+1}\beta(4j+1)$$
$$={1 \over 4}\left\{\cot\left[{1 \over 4}(\pi+(1+i)\sqrt{2\pi y})\right]+
\tan\left[{1 \over 4}(\pi-(1-i)\sqrt{2\pi y})\right] \right.$$
$$\left. +\tan\left[{1 \over 4}(\pi+(1-i)\sqrt{2\pi y})\right]+
\tan\left[{1 \over 4}(\pi+(1+i)\sqrt{2\pi y})\right]\right\}$$
$$={1 \over 2}\left\{\mbox{sec}\left[(1+i)\sqrt{{\pi y} \over 2}\right]
+ \mbox{sech}\left[(1+i)\sqrt{{\pi y} \over 2}\right]\right\}.$$

\bigskip
{\bf 5.  Connection with the Stieltjes coefficients}
\bigskip

{\bf Proposition 10}.  (a)
$$\int_0^\infty {{\eta(ix)} \over x}dx=4 \mbox{coth}^{-1} \sqrt{3}.$$
(b) 
$$-\int_0^\infty \ln x{{\eta(ix)} \over x}dx=
{1 \over {\sqrt{3}}}\left[
2\sqrt{3}\ln(2-\sqrt{3})(\gamma+\ln(4\pi))-\gamma_1
\left({1 \over {12}}\right)+\gamma_1\left({5 \over {12}}\right)+\gamma_1
\left({7 \over {12}}\right)-\gamma_1\left({{11} \over {12}}\right)\right].$$

(a) We write (7) in \cite{glasser2009} in terms of the $L$-series
$$L_{12}(s)={1 \over {12^s}}\left[\zeta\left(s,{1 \over {12}}\right)-\zeta\left(
s,{5 \over {12}}\right)-\zeta\left(s,{7 \over {12}}\right)+\zeta\left(
s,{{11} \over {12}}\right)\right],$$
so that
$$\int_0^ \infty x^{-s}\eta(ix)dx={{8\sqrt{3}\pi} \over {(4\pi)^s}}
{{\Gamma(2s-1)} \over {\Gamma(s)}} L_{12}(2s-1). \eqno(5.1)$$
Taking the limit as $s \to 1$ in (5.1),
$$\int_0^ \infty x^{-1}\eta(ix)dx={1 \over {2\sqrt{3}}}
\left[-\psi \left({1 \over {12}}\right)+\psi\left({5 \over {12}}\right)+\psi\left({7 \over {12}}\right)-\psi\left({{11} \over {12}}\right)\right].$$
Upon simplication of the differences of the pairs of digamma function values,
we obtain the stated result.

(b)  We take the derivative of both sides of (5.1) and then put $s \to 1$.
Let $\gamma_1(a)$ be the first Stieltjes constant as appears in the regular
part of the Laurent expansion of the Hurwitz zeta function 
\cite{coffeyfirst,coffeyrmjm2011},
$$\zeta(s,a)={1 \over {s-1}}+\sum_{n=0}^\infty {{(-1)^n} \over {n!}}\gamma_n(a)(s-1)^n, ~~~~~~s \neq 1.$$
Then 
$$-\int_0^\infty \ln x{{\eta(ix)} \over x}dx = {1 \over {\sqrt{3}}}\left[
2\sqrt{3}\ln(2-\sqrt{3})(\gamma+\ln(576\pi))-\gamma_1
\left({1 \over {12}}\right)+\gamma_1\left({5 \over {12}}\right)+\gamma_1
\left({7 \over {12}}\right)-\gamma_1\left({{11} \over {12}}\right)\right].$$
\qed

{\it Remark}.  We could further use an evaluation based upon \cite{coffeyrmjm2011}
(Proposition 1)    
in order to replace the pairs of differences of $\gamma_1$ values with a linear
combination of values $\zeta''(0,k/12)$ with $k=1,5,7,11$.

The following may be useful in conjunction with the quasi-periodicity of $\eta$
shown in (1.1), $\eta(ix\pm j)=e^{\pm i \pi j/12}\eta(ix)$ for $j \in \mathbb{N}$.
\newline{\bf Lemma}. (a)
$$\int_0^1 \eta(ix)dx={{2 \pi} \over \sqrt{3}}-{{12} \over \pi}\sum_{n=-\infty}
^\infty {{(-1)^n} \over {(6n-1)^2}}e^{-(6n-1)^2 \pi/12}$$
and (b)
$$\int_0^1 \eta^3(ix)dx=1-{4 \over \pi}\sum_{n=0}^\infty {{(-1)^n} \over {(2n+1)}} e^{-(2n+1)^2\pi/4}.$$

These are easy consequences of (1.2) and (1.3), using the value $\beta(1)=\pi/4$. \qed

In addition, 
$$\int_j^{j+1}\eta(ix)dx={{12} \over \pi}\sum_{n=-\infty}^\infty \left[
e^{-j\pi(6n-1)^2/12}-e^{-(j+1)\pi(6n-1)^2/12}\right]{{(-1)^n} \over {(6n-1)^2}}.$$
Then
$$\sum_{j=0}^\infty \int_j^{j+1}\eta(ix)dx=\int_0^\infty \eta(ix)dx
={{12} \over \pi} \sum_{n=-\infty}^\infty {{(-1)^n} \over {(6n-1)^2}}
={{2\pi} \over \sqrt{3}}.$$
This method is quite distinct from how the evaluation (9) was obtained in
\cite{glasser2009}.

\bigskip
{\bf 6.  A future direction:  an example of another lacunary case of $\eta$}
\bigskip

A $q$ series is lacunary if the arithmetic density of its coefficients is zero.
There is a result of Serre \cite{serre} showing that the only even values of $d$ for
which $\eta^d(\tau)$ is lacunary are $d=2,4,6,8,10,14$, and $26$.  It seems to be
unknown whether there are any odd values of $d$, other than $1$ and $3$, for which
$\eta^d(\tau)$ is lacunary. 

The following holds in the case of $d=14$.
\newline{\bf Proposition 11}.
$$-180\sqrt{3}\int_0^\infty \eta^{14}(ix)dx
={1 \over {24\pi}}\sum_{m_1,m_2=0}^\infty {{(-1)^{m_1} \mbox{Im} [(6m_1+2+i\sqrt{3}(4m_2+1)^6]} \over {(6m_1+2)^2+3(4m_2+1)^2}}.$$

Based upon a mention by Winquist \cite{winquist},
$$b(a^2-b^2)a(a^2-9b^2)=a^5b-10a^3b^3+9ab^5={1 \over {6\sqrt{3}}}
\mbox{Im}[(a+ib\sqrt{3})^6],$$
giving
$$\sum_{a \equiv 2, mod 6;b\equiv 1, mod 4}(-1)^{(a-2)/6}\mbox{Im}
[(a+ib\sqrt{3})^6]q^{(a^2+3b^2)/12}=-180\sqrt{3}\eta^{14}(\tau).$$
Integrating both sides of this relation,
$$-180\sqrt{3}\int_0^\infty \eta^{14}(ix)dx={1 \over {24\pi}}\sum_{a \equiv 2, mod 6;b\equiv 1, mod 4}(-1)^{(a-2)/6}{{\mbox{Im} [(a+ib\sqrt{3})^6]} \over {(a^2+3b^2)}},$$
the result follows.  \qed




\bigskip
\centerline{\bf Acknowledgement}
Useful correspondence with M.\ L.\ Glasser is gratefully acknowledged.

\pagebreak


\begin{thebibliography}{99}
\bibitem{nbs}M. Abramowitz and I. A. Stegun,
{Handbook of Mathematical Functions, Washington, National Bureau of Standards (1964).}
\bibitem{andrews}G. E. Andrews,
{The theory of partitions, Cambridge University Press (1998).}
\bibitem{a2r}G. E. Andrews, R. Askey, and R. Roy,
{Special functions, Cambridge University Press (1999).}
\bibitem{apostolmodfcns}T. Apostol,
{Modular functions and Dirichlet series in number theory, Springer (1990).}
\bibitem{ramapart3}B. C. Berndt,
{Ramanujan's notebooks, part III, Springer (1990).}
\bibitem{coffeyfirst}M. W. Coffey,
{Series representations of the Riemann and Hurwitz zeta functions and series and integral representations of the first Stieltjes constant, arXiv:1106.5147 (2011).}
\bibitem{coffeyrmjm2011}M. W. Coffey
{On representations and differences of Stieltjes coefficients, and other relations,
Rocky Mtn. J. Math. {\bf 41}, 1815-1846 (2011).}
\bibitem{summatory}M. W. Coffey,
{Summatory relations and prime products for the Stieltjes constants, and other
results, arXiv:1701.07064 (2017).}
\bibitem{coffeyrama}M. W. Coffey,
{Functional relations for the Stieltjes constants, Ramanujan J. {\bf 39}, 577-601 (2016).}
\bibitem{fuchs}J. Fuchs, 
{Affine Lie algebras and quantum groups, Cambridge University Press (1992).}
\bibitem{glasser2009}M. L. Glasser,
{Some integrals of the Dedekind $\eta$-function, J. Math. Analy. Appls. {\bf 354},
490 (2009).}
\bibitem{grad}I. S. Gradshteyn and I. M. Ryzhik,
{Table of Integrals, Series, and Products, Academic Press, New York (1980).}
\bibitem{kac}V. Kac,
{Infinite-dimensional algebras, Dedekind's $\eta$-function, classical M\"{o}bius
function and the very strange formula, Adv. Math. {\bf 30}, 85-136 (1978).}
\bibitem{patkowski}A. E. Patkowski,
{Some remarks on Glaisher-Ramanujan type integrals, arXiv1505.01530v1 (2015);
A. E. Patkowski and M. Wolf, Comp. Methods Sci. Tech. {\bf 22}, 103-108 (2016).}
\bibitem{rroy2017}R. Roy,
{Elliptic and modular functions from Gauss to Dedekind to Hecke, Cambridge
University Press (2017).}
\bibitem{serre}J.-P. Serre,
{Sur la lacunarit\'{e} des puissances de $\eta$, Glasgow Math. J. {\bf 27}, 203-221 (1985).}
\bibitem{winquist}L. Winquist,
{An elementary proof of $p(11m+6) \equiv 0$ (mod $11$), J. Combin. Theory {\bf 6},
56-59 (1969).}
\end{thebibliography}
\end{document}